\documentclass[11pt, a4paper]{article}

\usepackage{amsmath, amsthm, amssymb, amsfonts}

\usepackage[english]{babel}
\usepackage[Glenn]{fncychap}
\usepackage{comment}
\usepackage{graphicx}
\usepackage{fancyhdr}
\usepackage{scrextend}
\usepackage[figurename=Figure]{caption}

\usepackage{color}

\newtheorem{defi}{Definition}
\newtheorem{lem}{Lemma}
\newtheorem{thm}{Theorem}

\begin{document}
\begin{center}
{\Large\bf  Monotone-iterative technique for an initial value problem for difference equations  with non--instantaneous impulses}
\end{center}
\vspace {1 cm} 

\begin{center} { S. HRISTOVA  $^{a}$} \end{center}

\vspace {1 cm}

\begin{quote}
 $ ^{a}$ \  University of Plovdiv, \ 24 Tzar Asen 24, \ 4000 Plovdiv, \ Bulgaria
 \\ {\it   e-mail  address}:\ snehri@uni-plovdiv.bg 

\end{quote}

\bigskip

\begin{abstract}
In this paper a special type of difference equations is  investigated. The impulses start abruptly at some points  and their  action continue on given finite  intervals. This
type of equations is used to model a real process. An algorithm,
namely, the monotone iterative technique is suggested to solve the initial  value problem for nonlinear difference equations with non-instantaneous impulses 
approximately. An important feature of our algorithm is that each successive
approximation of the unknown solution is equal to the  unique solution of an appropriately
constructed initial value problem for a linear difference equation with with non-instantaneous impulses,
and a formula for its explicit form is given.  It is proved  both  sequences are convergent  and their limits are minimal and maximal solutions  of the considered problem.

\end{abstract}


\section{Introduction}

 In the real world life  there are many processes and phenomena that
are characterized by rapid changes in their state. In the literature
there are two popular types of impulses:
\begin{itemize}
\item[-]  {\it instantaneous impulses}- the duration of these changes is relatively short compared to the
overall duration of the whole process. The model is given by
impulsive differential equations (see, for example,  the monographs
\cite{S}, \cite {LB}, \cite{SP} and the cited references therein);
\item[-] {\it noninstantaneous impulses} -  an impulsive action, which starts at an
arbitrary fixed point and remains active on a finite time interval.
E. Hernandez and  D. O'Regan  (\cite {HR}) introduced this new class
of abstract differential equations where the impulses are not
instantaneous and they investigated the existence of mild and
classical solutions.
\end{itemize}

One of the problems in difference equations, which unknown function is involved in the present time on both parts of the equation nonlinearly, is the the obtaining of the solution. Often it could be done in a closed frm but approximately.One of the approximate method is based on  the method of upper and lower solutions, combined with a monotone-iterative technique. It is used to construct two monotonous sequences of upper and lower solutions of the nonlinear non-instantaneous impulsive difference equation. This method is applied for for difference equations in \cite {PA}, \cite {P} and for impulsive difference equations in \cite {AG}.

The idea of the method is to use upper and lower solutions for initial iteration and to construct successive approximation from the corresponding non-instantaneous impulsive linear equation. These functional sequences converge monotonically to the minimal and maximal solutions of the nonlinear equation.

\section{Statement of the problem}

Let  $\mathbb{Z}_+$ denote the set of all nonegative integers. Let the increasing sequence  $\{n_i\}_{i=0}^{p+1}:\  n_i\in\mathbb{Z}_+,\ n_i\geq n_{i-1}+3, i=1,2,\dots,p$ and the sequence  $ \{d_i\}_{i=1}^{p}:\ d_i\in\mathbb{Z}_+,\ 1 \leq d_i\leq n_{i+1}-n_i-2,\ i=1,2,\dots,p$ be given.
We denote  $\mathbb{Z}[a,b]= \{ z \in \mathbb{Z}_+ : a \leq z \leq b \},\ a, b \in \mathbb{Z}_+,\ a<b$ and $I_{k}=\mathbb{Z}[n_k+d_k,n_{k+1}-2],\ k \in \mathbb{Z}[0,p-1],\ I_p=\mathbb{Z}[n_p+d_p,n_{p+1}-1]$ and $J_k=\mathbb{Z}[n_k+1,n_k+d_k],\ k \in \mathbb{Z}[1,p]$  where $d_0=0$.

Consider the {\it initial value problem (IVP)} for the nonlinear {\it non--instantaneous impulsive difference equation (NIDE)}
\begin{equation}
\begin{split}
\label{EMeq1}
&x(n+1) =f(n,x(n),x(n+1)) \ \mbox {for} \ n\in \bigcup_{k=0}^{p}I_{k}, \\
&x(n_k) =F(k,x(n_k-1)),\ k \in \mathbb{Z}[1,p]\\
&x(n) =g(n,x(n),x(n_k)) \ \mbox {for} \ n\in  \bigcup_{k=1}^{p}J_{k},\\
&x(n_0) =x_0, 
\end{split}
\end{equation}        
where 
$x,x_0 \in \mathbb{R}$, $f : \bigcup_{k=0}^{p}I_{k} \times \mathbb{R} \times \mathbb{R} \to \mathbb{R},\ F : \mathbb{Z}[1,p]\times \mathbb{R} \to \mathbb{R},$ and $g : \bigcup_{k=1}^{p}J_{k} \times \mathbb{R} \times \mathbb{R}\to \mathbb{R}.$
\section{Preliminaries results}
\begin{defi}
We will say that the function $\alpha : \mathbb{Z}[n_0,n_{p+1}] \rightarrow \mathbb{R}$ is a minimal(maximal) solution of the IVP for NIDE \eqref{EMeq1} in $\mathbb{Z}[n_0,n_{p+1}]$ if it is a solution of \eqref{EMeq1} and for any solution $u(n), n \in \mathbb{Z}[n_0,n_{p+1}]$ the inequality $\alpha(n) \leq u(n)$ $(\alpha(n) \geq u(n))$ holds on $\mathbb{Z}[n_0,n_{p+1}]$ 
\end{defi}

                                                                                    \begin{defi}
The function $\alpha: \mathbb{Z}[n_0,n_{p+1}] \rightarrow \mathbb{R}$ is called lower (upper) solutions of IVP for NIDE \eqref{EMeq1}, if:
\begin{displaymath}
\begin{split}
&\alpha (n+1) \leq (\geq)f(n,\alpha(n),\alpha(n+1)), \ \mbox {for} \ n\in \bigcup_{k=0}^{p}I_{k},\\
&\alpha(n_k) \leq (\geq)F(k,\alpha(n_k-1)),\ k \in \mathbb{Z}[1,p]\\
&\alpha(n) \leq (\geq) g(n,\alpha(n),\alpha(n_k)), \ \mbox {for} \ n\in \bigcup_{k=1}^{p}J_{k}\\
&\alpha(n_0) \leq (\geq) x_0
\end{split}
\end{displaymath}
\end{defi} 
Consider the linear NIDE of the type
\begin{equation}
\begin{split}
\label{EMeq2}
& u(n+1) = Q_n u(n) + \sigma_n,\ \ n \in \bigcup_{k=0}^{p}I_{k},\\
&u(n_k) = T_k u(n_k-1) + \mu_k,\  k \in \mathbb{Z}[1,p]\\
&u(n) = M_n u(n) + L_n u(n_k) + \gamma_n, \ \ n \in \bigcup_{k=1}^{p}J_{k}\\
&u(n_0)= x_0,
\end{split}
\end{equation}
where $u,x_0 \in \mathbb{R}$, $Q_n:\ n \in \bigcup_{k=0}^{p}I_{k}$, $\sigma_n:\ n \in \bigcup_{k=0}^{p}I_{k}$,  $ L_n,\ M_n\not =1,\gamma_n:\ n\in\bigcup_{k=1}^{p}J_{k}$, and $T_k,\mu_k:\ k \in \mathbb{Z}[1,p]$ are given real constants.
\begin{lem}\label{EMlem1} The IVP for NIDE \eqref{EMeq2} has an unique solution  given by
\begin{equation}
\begin{split}
\label{EMsol2}
u(n) &= N(n)\sum_{j=n_0-1}^{n-1}(\prod_{i=j+1}^{n} R(i))\sigma_{j}\prod_{i=j+1}^{n-1}Q_i+N(n)\sum_{j=n_0}^{n}(\prod_{i=j+1} ^{n}R(i))\zeta(j)\prod_{i=j}^{n-1}Q_i\\
&+\tau(n)\ \ \mbox{for}\ n \in \mathbb{Z}[n_0,n_{p+1}]
\end{split}
\end{equation}
where $\sigma_{n_0-1}=x_0$,  $\sigma_n=0,\ Q_n=1$ for $n\in \mathbb{Z}[n_0,n_{p+1}]/\bigcup_{k=0}^{p}I_{k}$, 
\begin{equation}\label{EM16}\begin{split}
 N(n)=%
  \begin{cases}
    \frac{L_n}{1-M_n} \ \ \  \mbox {for}\ \  n \in \bigcup_{k=1}^{p}J_{k}\\
   1 \ \ \  \mbox {otherwise}
  \end{cases}
    \end{split}\end{equation}
		\begin{equation}\label{EM161}\begin{split}
 \tau(n)=%
  \begin{cases}
    \frac{\gamma_n}{1-M_n} \ \ \  \mbox {for}\ \  n \in \bigcup_{k=1}^{p}J_{k}\\
   0 \ \ \  \mbox {otherwise}
  \end{cases}
    \end{split}\end{equation}
		\begin{equation}\label{EM17}\begin{split}
 R(n)=%
  \begin{cases}
    N(n_k+d_k)=\frac{L_{n_k+d_k}}{1-M_{n_k+d_k}} \ \ \  \mbox {for}\ \  n = n_k+d_k+1,\ k \in \mathbb{Z}[1,p]\\
		T_k \ \ \  \mbox {for}\ \  n = n_k,\ \ \ \ \ k \in \mathbb{Z}[1,p]\\
   1 \ \ \  \mbox {otherwise}
  \end{cases}
    \end{split}\end{equation}
		\begin{equation}\label{EM171}\begin{split}
 \zeta(n)=%
  \begin{cases}
  \tau(n_k+d_k)= \frac{\gamma_{n_k+d_k}}{1-M_{n_k+d_k}} \ \ \  \mbox {for}\ \  n = n_k+d_k+1,\ k \in \mathbb{Z}[1,p]\\
		\mu_k \ \ \  \mbox {for}\ \  n = n_k,\ \ \ \ \ k \in \mathbb{Z}[1,p]\\
   0 \ \ \  \mbox {otherwise}
  \end{cases}
    \end{split}\end{equation}

\end{lem}
\textbf{Proof:} We will use an induction with respect to the interval. Let $n \in I_0=\mathbb{Z}[n_0,n_1-2].$ Then   we obtain
$u(n)= \sum_{j=n_0-1}^{n-1}\sigma_{j}\prod_{i=j+1}^{n-1}Q_i,\ n \in \mathbb{Z}[n_0 +1,n_1-1].$

Let $n=n_1.$ Then using $\sigma_{n_1-1}=0,\ \ Q_{n_1-1}=1$ we get
\begin{displaymath}
\begin{split}
u(n_1)&=T_1 u(n_1-1)+\mu_1=T_1 \sum_{j=n_0-1}^{n_1-1}\sigma_{j}\prod_{i=j+1}^{n_1-1}Q_i+\mu_1\\
&=  \sum_{j=n_0-1}^{n_1-1}(\prod_{i=j+1} ^{n_1}R(i))\sigma_{j}\prod_{i=j+1}^{n_1-1}Q_i+\mu_1.
\end{split} \end{displaymath}
Let $n \in J_1=\mathbb{Z}[n_1+1,n_1+d_1].$ Then using $\sigma_{j}=0,\ Q_{j}=1,$
$j \in \mathbb{Z}[n_1,n_1+d_1]$ we get
\begin{displaymath}
\begin{split}
u(n)&=\frac{L_n}{1-M_n}\sum_{j=n_0-1}^{n-1}(\prod_{i=j+1} ^{n_1}R(i))\sigma_{j}\prod_{i=j+1}^{n-1}Q_i+ \frac{L_n}{1-M_n}\mu_1+\frac{\gamma_n}{1-M_n}.\\
&=  N(n)\sum_{j=n_0-1}^{n-1}\sigma_{j}(\prod_{i=j+1}^{n}R(i))\prod_{i=j+1}^{n-1}Q_i+  N(n)\mu_1+\tau(n),\ n \in J_1.
\end{split} \end{displaymath}\normalfont
Let $n \in I_1=\mathbb{Z}[n_1+d_1,n_2-2].$ Then
\begin{displaymath}
\begin{split}
u(n)&= N(n_1+d_1)\sum_{j=n_0-1}^{n_1+d_1-1}\sigma_{j}(\prod_{i=j+1}^{n_1+d_1}R(i))\prod_{i=j+1}^{n_1+d_1-1}Q_i\prod_{i=n_1+d_1}^{n-1}Q_i\\
&\ \ \ +\mu_1N(n_1+d_1)\prod_{i=n_1+d_1}^{n-1}Q_i+\tau(n_1+d_1)\prod_{i=n_1+d_1}^{n-1}Q_i\\&\ \ \ +\sum_{j=n_1+d_1}^{n-1}\sigma_{j}\prod_{i=j+1}^{n-1}Q_i\\
&= \sum_{j=n_0-1}^{n-1}(\prod_{i=j+1}^{n}R(i))\sigma_{j}\prod_{i=j+1}^{n-1}Q_i+\sum_{j=n_0}^{n}(\prod_{i=j+1}^{n}R(i))\zeta(j)\prod_{i=j}^{n-1}Q_i.
\end{split}
\end{displaymath}
Let $n=n_2.$ Then we get 
\begin{displaymath}
\begin{split}
u(n_2)&=T_1 u(n_2-1)+\mu_2\\
&=  \sum_{j=n_0-1}^{n_2-1}(\prod_{i=j+1}^{n_2}R(i))\sigma_{j}\prod_{i=j+1}^{n_2-1}Q_i+\sum_{j=n_0}^{n_2}(\prod_{i=j+1}^{n_2}R(i))\zeta(j)\prod_{i=j}^{n_2-1}Q_i.
\end{split} \end{displaymath}
Let $n \in J_2=\mathbb{Z}[n_2+1,n_2+d_2].$ Then
\begin{displaymath}
\begin{split}
u(n)&=N(n)u(n_2)+\tau(n)\\
&=N(n)\sum_{j=n_0-1}^{n_2-1}(\prod_{i=j+1}^{n_2}R(i))\sigma_{j}\prod_{i=j+1}^{n_2-1}Q_i+N(n)\sum_{j=n_0}^{n_2}(\prod_{i=j+1}^{n_2}R(i))\zeta(j)\\&\ \ \ \times \prod_{i=j}^{n_2-1}Q_i+\tau(n)\\
&= N(n)\sum_{j=n_0-1}^{n-1}(\prod_{i=j+1}^{n}R(i))\sigma_{j}\prod_{i=j+1}^{n-1}Q_i+N(n)\sum_{j=n_0}^{n}(\prod_{i=j+1}^{n}R(i))\zeta(j)\\
&\ \ \ \times\prod_{i=j}^{n-1}Q_i+\tau(n).
\end{split}
\end{displaymath}
Continue this process step by step w.r.t. the interval by induction we proves the solution of NIDE \eqref{EMeq2} is given by \eqref{EMsol2} for all $n \in \mathbb{Z}[n_0,n_{p+1}].$ $\Box$
\begin{lem}\label{EMlem2}
Assume $m : \mathbb{Z}[n_0,n_{p+1}] \rightarrow \mathbb{R}$ satisfies the inequalities
\begin{equation}
\label{EMlem2i}
\begin{split}
&m(n+1) \leq Q_n m(n),\ \ n \in \bigcup_{k=0}^{p}I_{k},\\
&m(n_k) \leq T_k m(n_k-1), \ \ k \in \mathbb{Z}[1,p]\\
&m(n) \leq M_n m(n) + L_n m(n_k), \ \ n \in \bigcup_{k=1}^{p}J_{k}\\
&m(n_0) \leq 0,
\end{split}
\end{equation}
where $Q_n>0\Big{(} n \in \bigcup_{k=0}^{p}I_{k}\Big{)},$  $\ T_k>0\Big{(} k \in \mathbb{Z}[1,p]\Big{)}$ and $L_n > 0,\ M_n < 1\Big{(} n \in \bigcup_{k=0}^{p}J_{k}\Big{)}.$

Then $m(n) \leq 0$ for every $ n \in \mathbb{Z}[n_0,n_{p+1}].$
\end{lem}

The proof is based on an induction w.r.t. the interval and we
omit it.

\section{Main results}

For any pair of function $\alpha,\beta : \mathbb{Z}[n_0,n_{p+1}] \rightarrow \mathbb{R}$ such that $\alpha(n) \leq \beta(n)$ for $n \in \mathbb{Z}[n_0,n_{p+1}]$ we define the sets
\begin{displaymath}
\begin{split}
&S(\alpha,\beta)= \{ u : \mathbb{Z}[n_0,n_{p+1}] \rightarrow \mathbb{R}: \alpha(n) \leq u(n) \leq \beta(n), \ \ n \in \mathbb{Z}[n_0,n_{p+1}] \} \\
&\Omega_1(\alpha,\beta)= \{ u\in  \mathbb{R}: \alpha(n) \leq u \leq \beta(n), \ \ n \in  \bigcup_{k=0}^{p}I_{k}\}\\
&\Omega_2(\alpha,\beta)= \{ u \in \mathbb{R}: \alpha(n+1) \leq u \leq \beta(n+1), \ \ n \in  \bigcup_{k=0}^{p-1}I_{k}\}\\
&\Lambda(\alpha,\beta)= \{ u \in \mathbb{R}: \alpha(n) \leq u \leq \beta(n), \ \ n \in  \bigcup_{k=1}^{p}J_{k}\}\\
&\Gamma(\alpha,\beta)= \{ y \in \mathbb{R}: \alpha(n_k) \leq y \leq \beta(n_k), \ \ k \in \mathbb{Z}[1,p] \}\\
&\Upsilon(\alpha,\beta)= \{ z \in \mathbb{R}: \alpha(n_k-1) \leq z \leq \beta(n_k-1), \ \ k \in \mathbb{Z}[1,p] \}
\end{split}
\end{displaymath}

\begin{thm}
\label{EMT1}
Let the following conditions be fulfilled:
\begin{enumerate}
	\item  The functions $\alpha,\beta : \mathbb{Z}[n_0,n_{p+1}] \rightarrow \mathbb{R}$ are lower and upper solutions of the IVP for NIDE \eqref{EMeq1} and
$\alpha(n) \leq \beta(n)$ for $n \in \mathbb{Z}[n_0,n_{p+1}]$.	
	\item The function $f : \bigcup_{k=0}^{p}I_{k} \times \Omega_1(\alpha,\beta) \times \Omega_2(\alpha,\beta)$ 
is continuous in its second and third arguments  and there exist functions  $K: \bigcup_{k=0}^{p}I_{k} \to (-\infty, 1)$ and $P: \bigcup_{k=0}^{p}I_{k} \to (0,\infty)$ such that for any $n \in \bigcup_{k=0}^{p}I_{k}$ and $x_1,x_2 \in \Omega_1(\alpha,\beta),$ with $x_1 \leq  x_2,$ and $x_3,x_4 \in \Omega_2(\alpha,\beta),$ with $x_3 \leq  x_4$ the inequality
\begin{displaymath}
f(n,x_1,x_3) - f(n,x_2,x_4) \leq P(n)(x_1 - x_2) + K(n)(x_3 - x_4)
\end{displaymath}
holds.
	\item The function $F: \mathbb{Z}[1,p]\times\mathbb{R}\to \mathbb{R}$
is continuous in its second argument  and there exists a function $T: \ \mathbb{Z}[1,p]\to (0,\infty) $ such that for any $k \in \mathbb{Z}[1,p]$ and $z_1,z_2 \in \Upsilon(\alpha,\beta)$ with $z_1 \leq  z_2$ 
\begin{displaymath}
F(k,z_1) - F(k,z_2) \leq T(k)(z_1-z_2).
\end{displaymath}
	\item The function $g : \bigcup_{k=1}^{p}J_{k} \times \Lambda(\alpha,\beta) \times \Gamma(\alpha,\beta)\to\mathbb{R}$
is continuous in  its second and third arguments  and there exist  functions  $M: \bigcup_{k=1}^{p}J_{k}\to (-\infty, 1)$  and $L: \bigcup_{k=1}^{p}J_{k}\to (0,\infty)$  such that for any $n \in \bigcup_{k=1}^{p}J_{k}$ and $y_1,y_2 \in \Lambda(\alpha,\beta),$ with $y_1 \leq  y_2,$ and $y_3,y_4 \in \Gamma(\alpha,\beta),$ with $y_3 \leq  y_4$ the inequality
\begin{displaymath}
g(n,y_1,y_3) - g(n,y_2,y_4) \leq M(n)(y_1 - y_2) + L(n)(y_3 - y_4)
\end{displaymath}
holds.
\end{enumerate}

Then there exist two sequences of discrete functions $\{\alpha^{(j)}(n)\}_{0}^{\infty}$ and $\{\beta^{(j)}(n)\}_{0}^{\infty}$, $n\in \mathbb{Z}[n_0,n_{p+1}]$ with $\alpha^{(0)}=\alpha$ and $\beta^{(0)}=\beta$ such that:

a) The sequences are nondecreasing and nonincreasing, respectively and
\begin{displaymath}
\alpha(n) \leq \alpha^{(j)}(n) \leq \beta^{(j)}(n) \leq \beta(n),\ \textrm{for}\ \ n \in \mathbb{Z}[n_0,n_{p+1}];
\end{displaymath}

b) The functions $\alpha^{(j)}(n)$ and $\beta^{(j)}(n)$ are lower and upper solutions of the IVP for NIDE \eqref{EMeq1}, respectively;

c) Both sequences are convergent on $\mathbb{Z}[n_0,n_{p+1}]$;

d) The limits $\lim_{j \to \infty}\alpha^{(j)}(n)=A(n),$ $\lim_{j\to\infty}\beta^{(j)}(n)=B(n)$ are the minimal and maximal solutions of IVP for NIDE \eqref{EMeq1} in $S(\alpha,\beta),$ respectively;

e) If IVP for NIDE \eqref{EMeq1} has an unique solution $u(n) \in S(\alpha,\beta)$, then $A(n)\equiv u(n) \equiv B(n)$ for $n \in \mathbb{Z}[n_0,n_{p+1}].$ 

\end{thm}

\textbf{Proof:} For any arbitrary fixed function $\eta \in S(\alpha,\beta),$ we consider the IVP for the linear NIDE 
\begin{equation}
\begin{split}
\label{EMeq3}
& u(n+1) = P(n)u(n) + K(n) u(n+1) + \psi(n,\eta(n),\eta(n+1)),\ \ n \in \bigcup_{k=0}^{p}I_{k} \\
&u(n_k) = T(k) u(n_k-1) + \upsilon(k,\eta(n_k-1)), \ \ k \in \mathbb{Z}[1,p]\\
&u(n) = M(n) u(n) + L(n) u(n_k) + \xi(n,\eta(n),\eta(n_k)), \ \ n \in \bigcup_{k=1}^{p}J_{k}\\
&u(n_0)= x_0,
\end{split}
\end{equation}
where $u,x_0 \in \mathbb{R},$ and
\begin{equation}
\begin{split}
\label{EM10}
&\psi(n,x,y)=f(n,x,y)-P(n) x-K(n) y, \ \ n \in \bigcup_{k=0}^{p}I_{k}, \\
&\upsilon(k, x)=F(k,x)-T(k) x,\ k \in \mathbb{Z}[1,p]\\
&\xi(n,x,y)=g(n,x,y)-M(n) x - L(n)y, \ \ n \in \bigcup_{k=1}^{p}J_{k}.
\end{split}
\end{equation}

According to Lemma \ref{EMlem1} the IVP for linear NIDE \eqref{EMeq3} has an unique solution given by \eqref{EMsol2} with $\sigma_n=\frac{\psi(n,\eta(n),\eta(n+1))}{1-K(n)},\ Q_n=\frac{P(n)}{1-K(n)}$,  $\gamma_n=\xi(n,\eta(n),\eta(n_j))$ , $\mu_k=\upsilon(k,\eta(n_k-1))$, $T_k=T(k), M_n=M(n),L_n=L(n)$.

For any function $\eta \in S(\alpha,\beta)$ we define the operator $Q:S(\alpha,\beta) \rightarrow S(\alpha,\beta)$ by $Q\eta=u,$ where $u$ is the unique solution of IVP for the linear NIDE \eqref{EMeq3} for the function $\eta.$ 
The operator $Q$ has the following properties:

(P1) $\alpha \leq Q\alpha$, $\beta \geq Q\beta$

(P2) $Q$ is a monotone nondecreasing operator in $S(\alpha,\beta).$

To prove (P1) set $Q\alpha=\alpha^{(1)},$ where $\alpha^{(1)}$ is the unique solution of \eqref{EMeq3} with $\eta = \alpha$ and let $m(n)=\alpha(n)-\alpha^{(1)}(n),\ n \in \mathbb{Z}[n_0,n_{p+1}]$

For any $n \in  \bigcup_{k=0}^{p}I_{k} $ we obtain the inequality
\begin{displaymath}
\begin{split}
m(n+1) &= \alpha (n+1)-P(n)\alpha^{(1)}(n) - K(n) \alpha^{(1)}(n+1)\\
&\ \ \ -\psi(n,\alpha(n),\alpha(n+1))\\
&\leq P(n)(\alpha(n) - \alpha^{(1)}(n))+K(n)(\alpha(n+1) - \alpha^{(1)}(n+1))\\
&= P(n)m(n) + K(n) m(n+1).
\end{split}
\end{displaymath}
Hence the inequality $m(n+1)\leq \frac{P(n)}{1-K(n)}\ m(n)$ holds for $n \in \bigcup_{k=0}^{p}I_{k} .$

For any $n=n_k$ we obtain
\begin{displaymath}
\begin{split}
m(n_k) &\leq F(k,\alpha(n_k-1)) - T(k)\alpha^{(1)}(n_k-1)- F(k,\alpha(n_k-1))\\&\ \ \ + T(k) \alpha(n_k-1)= T(k) m(n_k-1).
\end{split}
\end{displaymath}

For any $n \in \bigcup_{k=1}^{p}J_{k}$ we get
\begin{displaymath}
\begin{split}
m(n) &\leq g(n,\alpha(n),\alpha(n_k)) - M(n) \alpha^{(1)}(n) - L(n) \alpha^{(1)}(n_k)\\
&\ \ \ -g(n,\alpha(n),\alpha(n_k))+ M(n) \alpha(n) + L(n) \alpha(n_k)\\
&= M(n) m(n) + L(n) m(n_k)
\end{split}
\end{displaymath}

Therefore, the function $m(n)$ satisfies the inequalities \eqref{EMlem2i} with $Q_n=\frac{P(n)}{1-K(n)}$,  $T_k=T(k), M_n=M(n),L_n=L(n)$. According to Lemma \ref{EMlem2} the function $m(n)$ is non-positive in $\mathbb{Z}[n_0,n_{p+1}],$ i.e. $\alpha \leq Q\alpha.$
Analogously it can be proved that the inequality $\beta \geq Q \beta$ holds.

To prove (P2) we consider two arbitrary function $\eta_1,\ \eta_2 \in S(\alpha,\beta)$ such that $\eta_1(n) \leq \eta_2(n)$ for $n \in \mathbb{Z}[n_0,n_{p+1}].$
Let $u^{(1)}=Q\eta_1$  and $u^{(2)}=Q\eta_2.$ Denote $m(n)=u^{(1)}(n)-u^{(2)}(n),\\n \in \mathbb{Z}[n_0,n_{p+1}].$

For any $n \in \bigcup_{k=0}^{p}I_{k} $ we obtain the inequality
\begin{displaymath}
\begin{split}
m(n+1) &= P(n) u^{(1)}(n) + K(n) u^{(1)}(n+1) + f(n,\eta_1(n),\eta_1(n+1))\\ 
&\ \ \ - P(n) \eta_1(n)- K(n) \eta_1(n+1) - P(n) u^{(2)}(n) - K(n) u^{(2)}(n+1)\\
&\ \ \ - f(n,\eta_2(n),\eta_2(n+1))+ P(n) \eta_2(n) + K(n) \eta_2(n+1)\\
&\leq P(n) m(n) + K(n) m(n+1)
\end{split}
\end{displaymath}
Hence the inequality $m(n+1)\leq \frac{P(n)}{1-K(n)}\ m(n)$ holds for $n \in I_k,\ k \in \mathbb{Z}[0,p].$

For any $n=n_k,  k \in \mathbb{Z}[1,p]$ we get
\begin{displaymath}
\begin{split}
m(n_k) &=T(k,u^{(1)}(n_k-1)  - u^{(2)}(n_k-1)) - T(k)(\eta_1(n_k-1)-\eta_2(n_k-1))\\
&\ \ \ + F(k,\eta_1(n_k-1)) - F(k,\eta_2(n_k-1)) \leq T(k) m(n_k-1)
\end{split}
\end{displaymath}

For any $n \in \bigcup_{k=1}^{p}J_{k}$ we obtain
\begin{displaymath}
\begin{split}
m(n) &= M_(n) u^{(1)}(n) + L(n) u^{(1)}(n_k) + g(n,\eta_1(n),\eta_1(n_k))\\
&\ \ \ - M(n) \eta_1(n) - L(n) \eta_1(n_k) - M(n) u^{(2)}(n) - L(n) u^{(2)}(n_k)\\ 
&\ \ \ - g(n,\eta_2(n),\eta_2(n_k)) + M(n) \eta_2(n) + L(n) \eta_2(n_k)\\
&\leq M(n)m(n) + L(n)m(n_k)
\end{split}
\end{displaymath}

According to Lemma \ref{EMlem2} with $Q_n=\frac{P(n)}{1-K(n)},$ $T_k=T(k), M_n=M(n),$ $L_n=L(n)$ the function $m(n) \leq 0,$ i.e. $Q\eta_1 \leq Q\eta_2,$ for $\eta_1(n) \leq \eta_2(n), n \in \mathbb{Z}[n_0,n_{p+1}].$

Let $\eta \in S(\alpha,\beta)$ be a lower solution of \eqref{EMeq1}. We consider the function $Q\eta=m.$ According to the proved $\eta(n) \leq m(n), n \in \mathbb{Z}[n_0,n_{p+1}].$

For any $n \in \bigcup_{k=0}^{p}I_{k} $ we get the inequality 
\begin{equation}
\begin{split}
\label{EMineq4}
m(n+1) 
&=P(n) m(n) + K(n) m(n+1) + f(n,\eta(n),\eta(n+1))\\ &\ \ \ - P(n) \eta(n) - K(n) \eta(n+1)\\
&\leq f(n,m(n),m(n+1))
\end{split}
\end{equation}

For any $n=n_k, k \in \mathbb{Z}[1,p]$ we obtain
\begin{equation}
\begin{split}
\label{EMineq55}
m(n_k) & =F(k,m(n_k-1)) - F(k,m(n_k-1)) + T(k) m(n_k-1) \\
&\ \ \ + F(k,\eta(n_k-1))- T(k)\eta(n_k-1)\\
&\leq F(k,m(n_k-1))
\end{split}
\end{equation}

For any $n \in \bigcup_{k=1}^{p}J_{k}$ we obtain
\begin{equation}
\begin{split}
\label{EMineq5}
m(n) &= g(n,m(n),m(n_k)) - g(n,m(n),m(n_k)) + M(n) m(n)\\ &\ \ \ + L(n) m(n_k) + f(n,\eta(n),\eta(n_k))- M(n) \eta(n) - L(n) \eta(n_k)\\
&\leq g(n,m(n),m(n_k))
\end{split}
\end{equation}

Inequalities \eqref{EMineq4},\eqref{EMineq55} and \eqref{EMineq5} prove the function $m$ is a lower solution of NIDE \eqref{EMeq1}. Similarly, if $\eta \in S(\alpha,\beta)$ is an upper solution of NIDE \eqref{EMeq1} then the function $m=Q\eta$ is an upper solution of \eqref{EMeq1}.

We define the sequences of functions $\{\alpha^{(j)}(n)\}_{0}^{\infty}$ and $\{\beta^{(j)}(n)\}_{0}^{\infty}$ by the equalities $\alpha^{(0)}=\alpha ,\ \beta^{(0)}=\beta,\ \alpha^{(j)}=Q\alpha^{(j-1)},\ \beta^{(j)}=Q\beta^{(j-1)}$.
The functions $\alpha^{(s)}(n)$ and $\beta^{(s)}(n)$ satisfy the initial value problem \eqref{EMeq3} with $\eta(n)=\alpha^{(s-1)}(n)$ and $\eta(n)=\beta^{(s-1)}(n), \ n \in \mathbb{Z}[n_0,n_{p+1}]$, respecticely.

According to Lemma \ref{EMlem2} the following representations are valid:
\begin{equation}
\begin{split}
\label{EMsolll2}
&\alpha^{(s)}(n) = N(n)\sum_{j=n_0-1}^{n-1}(\prod_{i=j+1}^{n} R(i))\frac{\psi(j,\alpha^{(s-1)}(j),\alpha^{(s=1)}(j+1))}{1-K(j)}\\&\times\prod_{i=j+1}^{n-1}\frac{P(i)}{1-K(i)}+\sum_{j=n_0}^{n}(\prod_{i=j+1} ^{n}R(i))\zeta(j)\prod_{i=j}^{n-1}\frac{P(i)}{1-K(i)}+\tau(n),\\ &\mbox{for}\ n \in \mathbb{Z}[n_0,n_{p+1}],
\end{split}
\end{equation}
where $\tau(n)$ is given by \eqref{EM161} for $\gamma_n=\xi(n,\alpha^{(s-1)}(n),\alpha^{(s-1)}(n_j)),\ n \in \bigcup_{k=1}^{p}J_{k},$ $j \in \mathbb{Z}[1,p]$ and $\zeta(n)$ is given by \eqref{EM171} for $\gamma_{n_k+d_k}=\xi(n_k+d_k,\alpha^{(s-1)}(n_k+d_k),\alpha^{(s-1)}(n_k)),\ \mu_k=\upsilon(k, \alpha^{(s-1)}(n_k-1)),\ k \in \mathbb{Z}[1,p]$.
\begin{equation}
\begin{split}
\label{EMsolll3}
&\beta^{(s)}(n) = N(n)\sum_{j=n_0-1}^{n-1}(\prod_{i=j+1}^{n} R(i))\frac{\psi(j,\beta^{(s-1)}(j),\beta^{(s-1)}(j+1))}{1-K(j)}\\ &\times\prod_{i=j+1}^{n-1}\frac{P(i)}{1-K(i)}+\sum_{j=n_0}^{n}(\prod_{i=j+1} ^{n}R(i))\zeta(j)\prod_{i=j}^{n-1}\frac{P(i)}{1-K(i)}+\tau(n),\\ &\mbox{for}\ n \in \mathbb{Z}[n_0,n_{p+1}],
\end{split}
\end{equation}
where $\tau(n)$ is given by \eqref{EM161} for $\gamma_n=\xi(n,\beta^{(s-1)}(n),\beta^{(s-1)}(n_j)),$ $j \in \mathbb{Z}[1,p]$ and $\zeta(n)$ is given by \eqref{EM171} for $\gamma_{n_k+d_k}=\xi(n_k+d_k,\beta^{(s-1)}(n_k+d_k),\beta^{(s-1)}(n_k)),\ \mu_k=\upsilon(k,\beta^{(s-1)}(n_k-1)),\ k \in \mathbb{Z}[1,p]$.

According to the above proved, functions $\alpha^{(s)} (n)$ and $\beta^{(s)}(n)$ are lower and upper solutions of NIDE \eqref{EMeq1}, respectively and they satisfy for $n \in \mathbb{Z}[n_0,n_{p+1}] $ the following inequalities

\begin{equation}
\label{EMmon}
\alpha^{(0)}(n) \leq \alpha^{(1)}(n) \leq \ldots \leq \alpha^{(s)}(n) \leq \beta^{(s)}(n) \leq \ldots \leq \beta^{(1)} (n) \leq \beta^{(0)}(n)
\end{equation}

Both sequences of discrete functions being monotonic and bounded are convergent on $\mathbb{Z}[n_0, n_{p+1}].$

Let $A(n) = \lim_{s \to \infty} \alpha^{(s)}(n), B(n) = \lim_{s \to \infty} \beta^{(s)}(n).$

Take a limit in \eqref{EMsolll2} for $s \to \infty $ we obtain \eqref{EMsol2} with $u(n)=A(n),\ \sigma_j=\frac{\psi(j,A(j),A(j+1))}{1-K(j)},\ Q_i=\frac{P(i)}{1-K(i)},$
\begin{equation}
\begin{split}
\label{EMminsol}
&A(n) = N(n)\sum_{j=n_0-1}^{n-1}(\prod_{i=j+1}^{n} R(i))\frac{\psi(j,A(j),A(j+1))}{1-K(j)}\prod_{i=j+1}^{n-1}\frac{P(i)}{1-K(i)}\\ &+\sum_{j=n_0}^{n}(\prod_{i=j+1} ^{n}R(i))\zeta(j)\prod_{i=j}^{n-1}\frac{P(i)}{1-K(i)}+\tau(n),\ \mbox{for}\ n \in \mathbb{Z}[n_0,n_{p+1}],
\end{split}
\end{equation}
where $\tau(n)$ is given by \eqref{EM161} for $\gamma_n=\xi(n,A(n),A(n_j)),\ j \in \mathbb{Z}[1,p]$ and $\zeta(n)$ is given by \eqref{EM171} for $\gamma_{n_k+d_k}=\xi(n_k+d_k,A(n_k+d_k),A(n_k)),\ \mu_k=\upsilon(k, A(n_k-1)),\ k \in \mathbb{Z}[1,p]$.

From \eqref{EMminsol} it follows the function $A(n)$ is a solution of NIDE \eqref{EMeq1}.

Similarly, we prove the function $B(n)$ is a solution of NIDE \eqref{EMeq1}.

Let $u \in S(\alpha, \beta)$ be a solution of IVP for NIDE \eqref{EMeq1}.
From inequalities \eqref{EMmon} it follows there exists a natural number p such that $p \in \mathbb{N}:$
\begin{displaymath}
\alpha^{(p)}(n) \leq u(n) \leq \beta^{(p)}(n) \ \ \textrm{for} \ \ n \in \mathbb{Z}[n_0, n_{p+1}].
\end{displaymath}

We introduce the notation $m(n)=\alpha^{(p+1)}(n) - u(n), \ \ n \in \mathbb{Z}[n_0, n_{p+1}].$

For any $n \in \bigcup_{k=0}^{p}I_{k}$ we get the inequality
\begin{displaymath}
\begin{split}
m(n+1) &=P(n) \alpha^{(p+1)}(n) + K(n) \alpha^{(p+1)}(n+1) +f(n,\alpha^{(p)}(n),\alpha^{(p)}(n+1))\\ 
&\ \ \ - P(n) \alpha^{(p)}(n) - K(n) \alpha^{(p)}(n+1) - f(n,u(n),u(n+1))\\
&\leq  P(n) m(n) + K(n) m(n+1)
\end{split}
\end{displaymath}
Hence the inequality $m(n+1)\leq \frac{P(n)}{1-K(n)}\ m(n)$ holds for $n \in \bigcup_{k=0}^{p}I_{k}.$

For any $n=n_k,\ \ k \in \mathbb{Z}[1,p]$ we obtain
\begin{displaymath}
\begin{split}
&m(n_k) = T(k) \alpha^{(p+1)}(n_k-1) + T(k) u(n_k-1) - T(k) u(n_k-1) \\&\ \ \ + F(k,(\alpha^{(p)}(n_k-1))- T(k) \alpha^{(p)}(n_k-1) - F(k,u(n_k-1))\\
& \leq T(k) m(n_k-1)
\end{split}
\end{displaymath}

For any $n \in \bigcup_{k=1}^{p}J_{k}$ we obtain
\begin{displaymath}
\begin{split}
m(n) &= M(n) \alpha^{(p+1)}(n) + L(n) \alpha^{(p+1)}(n_k) + g(n,\alpha^{(p)}(n),\alpha^{(p)}(n_k))\\ 
&\ \ \ - M(n) \alpha^{(p)}(n) - K(n) \alpha^{(p)}(n_k) - g(n,u(n),u(n_k))\\
& \leq  M(n) m(n) + L(n) m(n_k)
\end{split}
\end{displaymath}

According to Lemma \ref{EMlem2} with $Q_n=\frac{P(n)}{1-K(n)},$ $T_k=T(k), M_n=M(n),L_n=L(n)$ the function $m(n)$ is nonpositive, i.e. $\alpha^{(p+1)}(n) \leq u(n), \ \ n \in \mathbb{Z}[n_0,n_{p+1}].$ Similarly $\beta^{(p+1)}(n) \geq u(n), \ \ n \in \mathbb{Z}[n_0,n_{p+1}],$ and hence $\alpha^{(j+1)} \leq u(n) \leq \beta^{(j+1)}, \ \ n \in \mathbb{Z}[n_0,n_{p+1}].$ Since $\alpha^{(0)}(n) \leq u(n) \leq \beta^{(0)}(n)$ this proves by induction that $\alpha^{(j)}(n) \leq u(n) \leq \beta^{(j)}(n), \ \ n \in \mathbb{Z}[n_0,n_{p+1}], $for every $j.$

Taking the limit as $j \to \infty $ we conclude $A(n) \leq u(n) \leq B(n), \ \ n \in \mathbb{Z}[n_0,n_{p+1}].$ Hence $A(n)$ and $B(n)$ are minimal and maximal solutions of IVP for NIDE \eqref{EMeq1}, respectively.

Let the IVP for NIDE \eqref{EMeq1} has an unique solution $u(n) \in S(\alpha,\beta).$

Then from above it follows $A(n) \equiv u(n) \equiv B(n), \ \ n \in \mathbb{Z}[n_0,n_{p+1}].$ 
$\Box$

\section{Conclusions}

An algorithm for approximate solving an initial value problem for a nonlinear difference equations with non-instantaneous impulses is given and theoretically studied. It is based on the application of the method of lower and upper solutions. The suggested algoritm is appropriate for computerized and easy application to study discrete dynamical models.


\begin{thebibliography}{99}
\bibitem{EMRPA} R. P. Agarwal, \textit{Difference equations and inequalities}, Singapore, National University of Singapore, 2000.


\bibitem {AG}
R. P. Agarwal, S. Hristova,  A. Golev, K. Stefanova, Monotone-iterative method for mixed boundary value problems for generalized difference equations with “maxima”, {\it J. Appl. Math. Comput.}{\bf43}, 1, (2013) 213–-233. 


\bibitem{EMSE} S. Elaydi, \textit{An introduction to difference equations}, San Antonio, Dept. Math., Trinity University, 2005.

\bibitem{HR}
E. Hernandez, D. O'Regan, On a new class of abstract impulsive differential equations, {\it  Proc. Amer. Math. Soc.}, {\bf 141},
(2013), 1641--1649.


\bibitem {S}
 S. Hristova, {\it Qualitative inestigations and approximate methods for impulsive equations}, Nova Sci. Publ. Inc., New York, 2009.


\bibitem{LB}
V. Lakshmikantham, D.D. Bainov, P.S. Simeonov, {\it Theory of Impulsive Differential Equations}, World Scientific, Singapore, 1989.


\bibitem{P}
C. V. Pao, Monotone iterative methods for finite difference system of reaction-diffusion equations, {\it Numerische Math.}, {\bf 46}, 4, (1985) 571–-586.

\bibitem{PA}
P.Y.H. Pang, R.P. Agarwal
, Monotone iterative methods for a general class of discrete boundary value problems, {\it Comput. Math. Appl.}, {\bf 28},  1–3, (1994), 243--254.

\bibitem{SP}
A. M. Samoilenko, N. A. Perestyuk, {\it Impulsive differential equations}, World Scientific, Singapore, 1995.

\bibitem{WT}
P. Wang, Sh. Tian, Yonghong Wu, Monotone iterative method for first-order functional difference equations with nonlinear boundary value conditions, 
{\it Appl.  Math. Comput.}  {\bf 203},  1, (2008),  266–-272.




\end{thebibliography}
\end{document}